\documentclass[11pt]{article}

\pdfoutput=1 
\usepackage[pdftex]{graphicx}
\DeclareGraphicsExtensions{.pdf, .jpg}
\usepackage{amssymb,amsmath}

\def\R{\mathbf{R}}

\def\vol{\mathrm{vol}}

\newtheorem{thm}{Theorem}
\newtheorem{lem}{Lemma}
\newtheorem{cor}[thm]{Corollary}

\newtheorem{ques}{Question}

\newenvironment{pf}{\begin{trivlist}\item[]{\bf Proof:\ }}
{\mbox{}\hfill\rule{.08in}{.08in}\end{trivlist}}

\title{Semianalyticity of isoperimetric profiles}
\author{Renata Grimaldi, Stefano Nardulli and Pierre Pansu}
\begin{document}

\maketitle

\begin{quote}
{\small 
ABSTRACT: It is shown that, in dimensions $<8$, isoperimetric profiles of compact real analytic Riemannian manifolds are semi-analytic.
\par\medskip
RESUM\'E : On montre qu'en dimensions $<8$, le profil isop\'erim\'etrique d'une vari\'et\'e riemannienne compacte est semi-analytique.}
\end{quote}

\footnotetext[1]{Partially supported by Projet "Internazionalizzazione" "Propri\`et\`a asintotiche di variet\`a e di gruppi discreti" of Miur of Italy} 
\footnotetext[2]{P. Pansu, Univ Paris-Sud, Laboratoire de
Math\'ematiques d'Orsay, Orsay, F-91405} 
\footnotetext[3]{\hskip42pt CNRS, Orsay, F-91405.}

\section{Introduction}

\subsection{The problem}

Let $M$ be a compact real analytic Riemannian manifold. We are concerned with the regularity of the \emph{isoperimetric profile} of $M$.

Given $0<v<\vol(M)$, consider all integral currents in $M$ with volume $v$. Define $I_{M}(v)$ as the least upper bound of the boundary volumes of such currents. In this way, one gets a function $I_{M}:(0,\vol(M))\to\R_{+}$ called the \emph{isoperimetric profile} of $M$. In fact, for each $0<v<\vol(M)$, there exist currents in $M$ with volume $v$ and boundary volume $I_M (v)$. Such minimizing currents will be called {\em bubbles}, for short. 

Here is a typical example. For $a>0$, let $S$ denote the circle of length $2\pi$. Let $M=S \times S$. Then the isoperimetric profile of $M$ is easily computed to be
\begin{eqnarray*}
I_{M}(v)=\begin{cases}
\sqrt{4\pi v}& \text{for }0<v\leq 4\pi, \\
4\pi & \text{for }4\pi\leq v\leq 4\pi(\pi-1), \\
\sqrt{4\pi(4\pi^2 -v)}& \text{for }4\pi(\pi-1)\leq v<4\pi^2 .
\end{cases}
\end{eqnarray*}

\begin{center}
\includegraphics[width=3in]{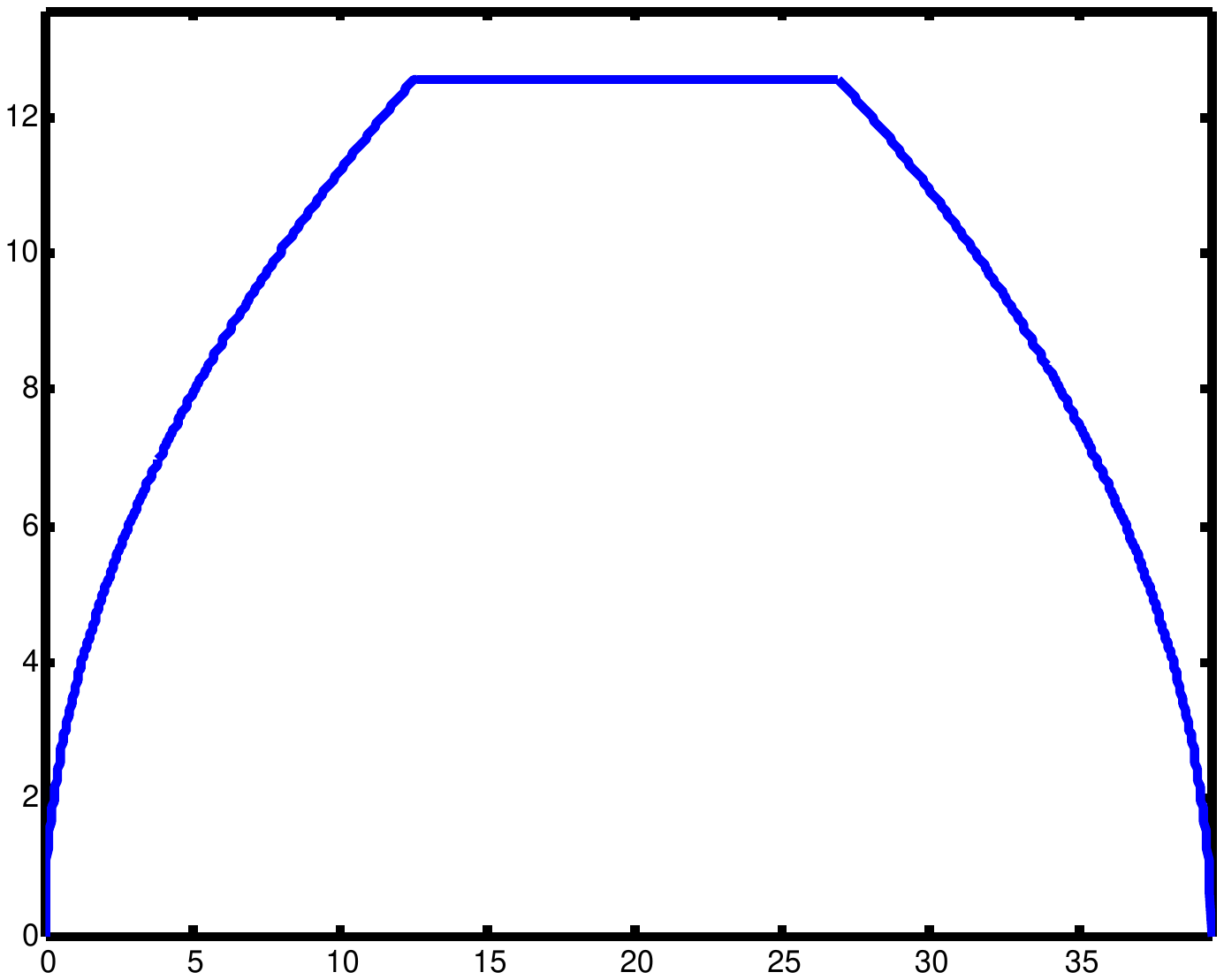}

\vskip -3cm Profile of a flat torus
\end{center}

This is proven as follows. In 2 dimensions, the boundaries of these bubbles are smooth, they have constant geodesic curvature, therefore they lift to disjoint unions of circles of equal radii or lines in $\R^2 =\tilde{M}$. It follows that bubbles are either round disks or annuli bounded by parallel geodesics, or complements of such. There remains to minimize boundary length among these three families.

\begin{ques}
For general real analytic manifolds, is it true that bubbles fall into finitely many analytic families, and that the profile is piecewise analytic ? 
\end{ques}

This has been proven in \cite{PP} in dimension 2 only.

\subsection{The results}

First, in a neighborhood of zero.

\begin{thm}
\label{small}
Let $M$ be a compact real analytic Riemannian manifold. There exists $\epsilon>0$ such that $I_M$ is real analytic on $(0,\epsilon)$.
\end{thm}

The isoperimetric profile of Euclidean space $\R^n$ is $I_{\R^n}(v)=n(\omega_n)^{1/n}v^{n-1/n}$, where $\omega_n$ is the volume of the unit ball in $\R^n$. In a curved manifold, $I_M (v)^{n-1}\sim n(\omega_n)^{1/n}v^{n-1/n}$ as $v$ tends to 0. 

\begin{ques}
For a compact analytic Riemannian $n$-manifold, is $I_M (v)$ an analytic function of $v^{1/n}$ on $[0,\epsilon)$ ?
\end{ques}

We have only a partial answer.

\begin{thm}
\label{racine}
Let $M$ be a compact real analytic Riemannian manifold. Assume that the absolute maxima of scalar curvature are nondegenerate critical points. Then there exists an analytic function $f$ defined in a neighborhood of 0 such that $I_M (v)=f(v^{1/n})$ for $v$ small enough.
\end{thm}

Away from 0, our result also requires an extra assumption.

\begin{thm}
\label{big}
Let $M$ be a compact real analytic Riemannian manifold. Let $0<v_0 <\vol(M)$. Assume that all bubbles of volume $v_0$ are smooth. Then $I_M$ is semi-analytic on a neighborhood of $v_0$. 
\end{thm}

Since bubbles are known to be smooth in dimensions $<8$, \cite{Almgren}, it follows that

\begin{cor}
\label{8}
If the dimension of $M$ is less than 8, $I_M$ is semianalytic on $[0,\vol(M)]$.
\end{cor}

\begin{ques}
Our method of proof relies on the regularity of bubbles. Can this be circumvented ?
\end{ques}

\section{Proof of Theorem \ref{small}}

It relies on results from \cite{Npvo}. There, it is shown that small bubbles is a subset of a smooth finite dimensional family of domains called {\em pseudo-balls}. We merely need show that if the metric is real analytic, pseudo-balls form a compact, finite dimensional real analytic set, on which the volume and boundary volume functions are real analytic.

\bigskip

Pseudo-balls are solutions of a differential equations which is weaker than constancy of mean curvature, but to which the implicit function theorem can be applied. Specificly, for $k\geq 0$, consider the bundle $\mathcal{F}^{k}\to M$ whose fiber at $p\in M$ consists in $C^{k,\alpha}$ functions on the unit sphere in the tangent space $T_p M$. There is a smooth map $\Phi:\R\times\mathcal{F}^2 \to\mathcal{F}^0$ with the following properties.
\begin{enumerate}
  \item Let $r>0$, $p\in M$ and $x\in \mathcal{F}^{2}_{p}$. If the graph, in polar coordinates, of $r(1+x)$ has constant mean curvature, then $\Phi(r,p,x)=0$.
  \item For all $p\in M$, $\Phi(0,p,0)=0$.
  \item The differential of $\Phi$ restricted to the fibers is an isomorphism.
\end{enumerate}

\begin{lem}
\label{analytic}
$\Phi:\R\times\mathcal{F}^2 \to\mathcal{F}^0$ is a real analytic map.
\end{lem}

\begin{pf}
$\Phi(r,p,x)$ has two components, $\Phi(r,p,x)=(\mathcal{A}(r,p,x),Q\circ \Psi(r,p,x))$. The fact that $\Psi$ is smooth is shown in \cite{Npvo}, Proposition 1.1. The map $\mathcal{A}$ is defined implicitely, by an equation whose smoothness is established in \cite{Npvo}, Lemma 2.1. When the metric is real analytic, the same proofs show that $\mathcal{A}$ and $\Psi$ are real analytic.
\end{pf}

From the implicit function theorem, it follows that for all $p\in M$ and for $r$ small enough, the equation $\Phi(r,p,x)=0$ has a unique small solution $x=x(r,p)$ which depends analytically on $(r,p)\in \R\times M$. Theorem 6 of \cite{Npvo} asserts that there exist $r_0 >0$ and $v_0 >0$ such that every bubble of volume less than $v_0$ coincides with the domain $\mathcal{N}^{+}(p,r)$ bounded by the graph, in polar coordinates, of the function $r(1+x(r,p))$, for some $r\leq r_0$ and some $p\in M$. Therefore, for $v<v_0$,
\begin{eqnarray*}
I_M (v)=\min\{\vol(\partial\mathcal{N}^{+}(p,r))\,|p\in M,\,0\leq r\leq r_0 ,\,\vol(\mathcal{N}^{+}(p,r))=v\}.
\end{eqnarray*}
Define the {\em lower contour} $c(A)$ of a subset $A\subset \R^2$ as the function $v\mapsto \inf\{w\in\R\,|\,(v,w)\in A\}$. Then the restriction of $I_M$ to $[0,v_0]$ coincides with the lower contour of the image of the real analytic map
\begin{eqnarray*}
\Omega:[0,r_0]\times M\to\R^2 ,\quad (r,p)\mapsto (\vol(\partial\mathcal{N}(p,r)),\vol(\partial\mathcal{N}^{+}(p,r))).
\end{eqnarray*}
Since $[0,r_0]\times M$ is semi-analytic and compact, $\Omega([0,r_0]\times M)$ is a compact subanalytic set. Its lower contour is a subanalytic function. According to S. Lojasiewicz (Theorem 6.1 in \cite{BM}), subanalytic functions on the real line are semi-analytic. Semi-analytic functions are piecewise analytic, thus there exists $\epsilon>0$ such that $I_M$ is analytic on $(0,\epsilon)$. This complete the proof of Theorem \ref{small}.

\section{Proof of Theorem \ref{racine}}

For $(p,\rho)\in M\times\R_+$, let $\beta(p,\rho)$ denote the pseudo-ball defined by
\begin{eqnarray*}
\beta(p,\vol(\mathcal{N}^{+}(p,r))^{1/n})=\mathcal{N}^{+}(p,r).
\end{eqnarray*}
Since $\vol(\mathcal{N}^{+}(p,r))^{1/n}\sim \omega_n^{1/n} r$ is a 1 to 1 analytic function of $r$, the notation is unambiguous. Let
\begin{eqnarray*}
f(p,\rho)=f(\vol(\partial\beta(p,\rho)))=f_\rho (p).
\end{eqnarray*}
Then $f$ is real analytic. Furthermore, among pseudo-balls of volume $v=\rho^{n}$, bubbles are characterized as minima of $f_\rho$. The following expansion
\begin{eqnarray*}
f_\rho (p)=c_n \rho^{n-1}(1-\frac{1}{2n(n+2)}\omega_n^{-2/n}Sc(p)\rho^{2}+\mathcal{O}(\rho^{4}))
\end{eqnarray*}
is computed in \cite{Npvo}, Lemma 3.6, compare \cite{Ye}. If the absolute maxima $p_1$,...,$p_k$ of the scalar curvature function $Sc$ are non degenerate critical points, then each of them deforms into a critical point $p_i (\rho)$ of $f_\rho$ that depends analytically on $\rho$. Therefore (Theorem 8 in \cite{Npvo}), 
\begin{eqnarray*}
I_M (\rho^n)=\min_{i=1,\ldots,k} f_\rho (p_i (\rho)).
\end{eqnarray*}
There exists $\epsilon>0$ and $i$ such that the minimum is equal to $f_\rho (p_i (\rho))$ for all $\rho\in[0,\epsilon)$. Indeed, otherwise, some function $f_\rho (p_i (\rho))-f_\rho (p_j (\rho))$ would change sign infinitely many times near 0, contradiction. Thus the right hand side is analytic on $[0,\epsilon)$. This completes the proof of Theorem \ref{racine}.

\section{Proof of Theorem \ref{big}}

We follow M. Tamm's strategy, \cite{Tamm}. We aim at including bubbles in a parametrized analytic variety. We shall first do this in a neighborhood of a smooth bubble $B$ with volume $v_0$. Our first candidate is the set of domains whose boundary is a graph in normal exponential coordinates to $\partial B$ and has constant mean curvature. To decide wether this set is a submanifold in some function space, let us examine the mean curvature operator and its linearization.

\subsection{Pseudo-bubbles}

Let $B$ be a smooth bubble with volume $v_0$. Let $H_{B}:C^{2,\alpha}(\partial B)\to C^{0,\alpha}(\partial B)$ denote the operator which to a function $u$ on $\partial B$ associates the mean curvature of the graph of $u$ in normal exponential coordinates to $\partial B$. In particular, $H_{B}(0)=H(\partial B)=h_{B}$ is the constant mean curvature of $\partial B$. Let $L_{B}:C^{2,\alpha}(\partial B)\to C^{0,\alpha}(\partial B)$ denote its linearization at 0 (sometimes called the {\em Jacobi operator}).

\begin{lem} For all $v\in C^{2,\alpha}(\partial B)$
\begin{equation}\label{LB}
           L_{B}(v)=-\Delta _{\partial B}v -(||II_{\partial B}||^2+Ric(\nu))v,
\end{equation}
where $\Delta _{\partial B}v=div(\nabla v)$ is the Laplace operator on $\partial B$ ( with negative spectrum when taken on the round sphere), $||II_{\partial B}||^2$ is the Hilbert-Schmidt squared norm ($tr(A^t A)$ for a square matrix $A$) of the second fundamental form of $\partial B$ and $Ric (\nu)$ is the Ricci curvature of the ambient manifold $M$ in the direction $\nu$ of the unit outward normal vector to $\partial B$ evaluated at a point of $\partial B$. 
\end{lem}

\begin{pf}
We recall here formula (6) of 3.3 of \cite{Nreg}
\begin{eqnarray}\label{h-norm1hyp}  
H(u)&=&-div_{(\mathbb{S}^{n-1},g_u )}(\frac{\overrightarrow{\nabla}_{g_u} u}{W_{u}})
-\frac{1}{W_{u} ^2}<\nabla_{\overrightarrow{\nabla}_{g_u} u}(\frac{u\overrightarrow{\nabla}_{g_u}u}{W_{u}}),\overrightarrow{\nabla}_{g_u} u>_{g_u}\\ 
\nonumber &+&\frac{u^2 }{W_{u} ^3}II_\theta ^{u} (\overrightarrow{\nabla}_{g_u} u,\overrightarrow{\nabla}_{g_u} u)\\ 
\nonumber &-&\frac{1}{W_{u}}H_\theta ^{u} (u)+\frac{1}{W_{u}}<\overrightarrow{\nabla}_{g_u} (\frac{1}{W_{u}}),u\frac{\overrightarrow{\nabla}_{g_u} u}{W_{u}})>_{g_u}.
\end{eqnarray}
Here $\theta$ denotes the gradient of the signed distance function to $\partial B$. For the meaning of the other terms involved in (\ref{h-norm1hyp}), see \cite{Nreg}. The operator $L_B $ satisfies then 
$H(tv)=H_{\partial B}+tL_B +\mathcal{O}(t^2)$. At first observe that only the first and fourth term of (\ref{h-norm1hyp}) contribute to $L_B$. Denoting by $U(r)=\nabla\theta$ the shape operator of the equidistant hypersurfaces to $\partial B$ at distance $r$, we have $U(r)=U_0+U_1 r+\cdots$,  hence by using the Riccati equation satisfied by $U$ (see \cite{Chav}) we can compute $U_1=-U_0 ^2- R$ where $R$ is the curvature tensor. Now taking traces we get $H_\theta ^{tv}=tr(-U)=H_{\partial B, ext}+(tr(U_0 ^2)+Ric)t+\cdots$ where $H_{\partial B, ext}$ is the outward mean curvature of the boundary of $B$. Finally (\ref{LB}) follows easily. 
\end{pf}

$L_{B}$ is a selfadjoint elliptic operator, which has a discrete spectrum. Let $K_B =\mathrm{kernel}(L_B)$, $m(B)=\mathrm{dim}(K_B)$. 

If $L_B$ were invertible (i.e. $m(B)=0$), the implicit function theorem would imply that nearby domains with constant mean curvature boundary come in one analytic family parametrized by the value of mean curvature. 

Unfortunately, $L_B$ is not always invertible. Therefore, instead of solving $H(u)=h$, $h\in \R$, we shall solve 
\begin{eqnarray*}
\Phi_{B}(u,h)=P_B (H_{B}(u)-h)=0,
\end{eqnarray*}
where $P_B$ is the orthogonal projection onto the $L^2$-orthogonal complement $K_{B}^{\bot}$ of $K_{B}$ in $C^{0,\alpha}(\partial B)$. Then 
\begin{eqnarray*}
\Phi_{B}:C^{2,\alpha}(\partial B)\times\R\to K_{B}^{\bot}
\end{eqnarray*}
is a real analytic map, whose linearization at 0 is $P_B \circ L_{B}$. By construction, it is onto. In fact, the restriction of $P_B \circ L_{B}$ to $K_{B}^{\bot}$ is an isomorphism (see for example \cite{Besse}, page 464). Note that $\Phi_{B}(0,h_{B})=0$. The following variant of the implicit function theorem provides us with an open neighborhood $U_{B}$ of $(0,h_{B})$ in $C^{2,\alpha}(\partial B)\times\R$ in which the solutions of $\Phi_{B}(u,h)=0$ form a real analytic submanifold. We shall call such solutions {\em $B$-pseudo-bubbles}.

\begin{lem}
\label{tfi}
Let $E$, $P$ and $F$ be real analytic Banach manifolds, let $e_0\in E$, $f_0\in B$, $p_0 \in P$ be such that $\Phi(e_0 ,p_0)=f_0$. Let $\Phi:E\times P\to B$ be a real analytic map. Assume that the differential $d\Phi$ of $\Phi$ at $(e_0 ,p_0)$ in the direction of $E$ has a finite dimensional kernel $K\subset T_{e_{0}} E$, which admits a closed complement $K^{\bot}$. Assume that the restriction of $d\Phi$ to $K^{\bot}$ is invertible. Then, in a neighborhood of $(e_0 ,p_0)$ the solutions of equation $\Phi(e,p)=f_0$ form a real analytic submanifold parametrized by a neighborhood of $(0,p_0)$ in $K\times P$.
\end{lem}

\begin{pf}
Apply the implicit function theorem to $\Psi:E\times P\to F\times K$ defined by $\Psi(e,p)=(\Phi(e,p),\pi_K (e))$ where $\pi_K$ is a local submersion onto $K$.
\end{pf}

\subsection{Compactness in $C^{2,\alpha}$-topology}

Let us define the $C^{2,\alpha}$-topology on the space of domains with smooth boundary as follows: as neighborhoods of a smooth domain $\beta$, take all domains $S$ whose boundary is the graph, in normal exponential coordinates, of a $C^{2,\alpha}$-small function on $\partial \beta$. Using a result from \cite{Nreg}, we show that on smooth bubbles with volume close to $v_0$, the topologies induced by the $C^{2,\alpha}$-topology on smooth domains and the flat topology on currents coincide. 

\begin{lem}
\label{epsilon}
Let $B$ be a bubble of volume $v_0$. For all $\delta>0$, there exist $\epsilon>0$ such that if $\beta$ is a bubble of volume $\in[v_0 -\epsilon,v_0 +\epsilon]$ with $\vol(\beta\Delta B)<\epsilon$, then there exists a smooth function $u$ on $\partial B$ with $\parallel u \parallel_{C^{2,\alpha}}<\delta$ such that $\partial\beta$ is the graph in normal exponential coordinates of $u$. Conversely, the graph of a $C^{2,\alpha}$-small function on $\partial B$ bounds a current which is close to $B$ volumewise.
\end{lem}

\begin{pf}
By contradiction. Otherwise, there exists a sequence $\beta_j$ of bubbles with $\vol(\beta_j)\to v_0$ and $\vol(\beta_j \Delta B)\to 0$ such that $\partial\beta_j$ is not the normal exponential graph of a $C^{2,\alpha}$-small function on $\partial B$. Theorem 1 of  \cite{Nreg} asserts that for $j$ large enough, $\partial\beta_j$ is the graph in normal exponential coordinates of a function $u_j$ on $\partial B$ whose $C^{2,\alpha}$-norm tends to zero, contradiction. The last statement is obvious.
\end{pf}

\begin{lem}
\label{smooth}
Assume that all volume $v_0$ bubbles in $M$ are smooth. Then there exists $\epsilon>0$ such that if $\beta$ is a bubble of volume $\in[v_0 -\epsilon,v_0 +\epsilon]$, then $\beta$ is smooth, and $(\beta,h_{\beta})$ belongs to the open set $U(B)$ for some volume $v_0$ bubble $B$.
\end{lem}

\begin{pf}
By contradiction. Otherwise, there exists a sequence $\beta_j$ of bubbles with $\vol(\beta_j)\to v_0$ that avoids all $U(B)$. By compactness of integral currents of bounded boundary volume, we can assume that $\beta_j$ converges in flat norm to some integral current $B$. By continuity of volume and semi-continuity of boundary volume, $B$ is a bubble of volume $v_0$. By assumption, $B$ is smooth. Theorem 1 of  \cite{Nreg} asserts that for $j$ large enough, $\partial\beta_j$ is the graph in normal exponential coordinates of a $C^{2,\alpha}$-small smooth function $u_j$ on $\partial B$. Therefore $\beta_j$ is smooth and $(\beta_j,h_{\beta_j})$ belongs to $U(B)$, contradiction.
\end{pf}

\begin{lem}
\label{N}
Assume that all volume $v_0$ bubbles in $M$ are smooth. Then there exists $\epsilon>0$ such that the set $\mathcal{B}$ of pseudo-bubbles with volumes $\in[v_0 -\epsilon,v_0 +\epsilon]$ is contained in a finite union of compact semi-analytic pieces of finite dimensional real analytic manifolds, on which the volume and boundary volume functions are real analytic.
\end{lem}

\begin{pf}
It was just proven that the set $\mathcal{B}$ of bubbles with volumes $\in[v_0 -\epsilon,v_0 +\epsilon]$ is compact in flat topology and the set $\mathcal{BH}=\{(B,h_B)\,|\,B\in \mathcal{B}\}$ is covered by the sets $U(B)$. According to Lemma \ref{epsilon}, it is compact in $C^{2,\alpha}$-topology as well. Therefore, $\mathcal{BH}$ can be covered with finitely many open sets $U(B_1)$,...,$U(B_N)$. The set $\Psi B_i$ of $B_i$-pseudo-bubbles in $U(B_i)$ is an analytic submanifold. There exist compact semi-analytic subsets $W_i \subset \Psi B_i$ which suffice to cover $\mathcal{BH}$.
\end{pf}

The proof of Theorem \ref{big} is completed in the same manner as the proof of Theorem \ref{small}. In dimensions less than 8, F. Almgren has shown that all bubbles are smooth. Therefore the profile is semi-analytic in a neighborhood of every point of the closed interval $[0,\vol(M)]$. It follows that it is semi-analytic on this interval. This proves Corollary \ref{8}.

\bibliographystyle{plain}
\bibliography{anpr}

\begin{thebibliography}{1}

\bibitem{Almgren}
Frederick~J. Almgren.
\newblock {\em Existence and regularity almost everywhere of solutions to
  elliptic variational problems with constraints}.
\newblock Number 165 in Memoirs. Amer. Math. Soc., Providence, R.I., 1976.

\bibitem{Besse}
Arthur~L. Besse.
\newblock {\em Einstein Manifolds}, volume~10 of {\em Ergebnisse der Math.
  Grenz.}
\newblock Springer Verlag, 1987.

\bibitem{BM}
Edward Bierstone and Pierre~D. Milman.
\newblock Semianalytic and subanalytic sets.
\newblock {\em Inst. Hautes \'Etudes Sci. Publ. Math.}, 67:5--42, 1988.

\bibitem{Chav}
Isaac Chavel.
\newblock {\em {Riemannian geometry. A modern introduction. 2nd ed.}}
\newblock Cambridge Studies in Advanced Mathematics 98. Cambridge: Cambridge
  University Press. xvi, 471~p., 2006.

\bibitem{Npvo}
Stefano Nardulli.
\newblock The isoperimetric profile of a compact riemannian manifold for small
  volumes.
\newblock 2007.

\bibitem{Nreg}
Stefano Nardulli.
\newblock Regularity of solutions of the isoperimetric problem that are close
  to a smooth manifold.
\newblock 2007.

\bibitem{PP}
Pierre Pansu.
\newblock Sur la r\'egularit\'e du profil isop\'erim\'etrique des surfaces
  riemanniennes compactes.
\newblock {\em Ann. Inst. Fourier}, 48:247--264, 1998.

\bibitem{Tamm}
Martin Tamm.
\newblock Subanalytic sets in the calculus of variations.
\newblock {\em Acta Math.}, 146:167--199, 1981.

\bibitem{Ye}
Rugang Ye.
\newblock Foliation by constant mean curvature spheres.
\newblock {\em Pacific J. Math.}, 147(2):381--396, 1991.

\end{thebibliography}

Keywords : Isoperimetric inequality, Riemann surface, calibration, Morse theory. 

Mathematics Subject Classification : 
53C20, 
49Q20, 
14P15, 
32B20. 

\vskip1cm
\noindent
Renata Grimaldi, Stefano Nardulli\\
Universit\`a di Palermo\\
Dipartimento di Metodi e Modelli Matematici\\
Facolta di Ingegneria\\
Vialle delle Scienze - 90128 Palermo (Italia)\\
\smallskip\noindent
{\tt grimaldi@unipa.it, nardulli@unipa.it}
\par\medskip\noindent
Pierre Pansu\\
Laboratoire de Math\'ematique d'Orsay\\
UMR 8628 du C.N.R.S.\\
B\^atiment 425\\
Universit\'e Paris-Sud 11 - 91405 Orsay (France)\\
\smallskip\noindent
{\tt\small Pierre.Pansu@math.u-psud.fr}\\
http://www.math.u-psud.fr/$\sim$pansu

\end{document}